\documentclass[12pt]{article}
\usepackage{amssymb,amsthm,amsmath,latexsym}
\usepackage{hyperref}
\newtheorem{thm}{Theorem}
\newtheorem{prop}[thm]{Proposition}

\newtheorem{cor}[thm]{Corollary}
\theoremstyle{remark}
\newtheorem{rem}[thm]{Remark}
\newcommand{\FF}{\mathbb{F}}

\DeclareMathOperator{\wt}{wt}

\begin{document}

\title{There is no $[21,5,14]$ code over $\FF_5$}

\author{
Makoto Araya\thanks{Department of Computer Science,
Shizuoka University,
Hamamatsu 432--8011, Japan.
email: araya@inf.shizuoka.ac.jp}
and 
Masaaki Harada\thanks{
Department of Mathematical Sciences,
Yamagata University,
Yamagata 990--8560, Japan, and
PRESTO, Japan Science and Technology Agency, Kawaguchi,
Saitama 332--0012, Japan. 
email: mharada@sci.kj.yamagata-u.ac.jp}
}

\maketitle

\begin{abstract}
In this note, we demonstrate that there is no
$[21,5,14]$ code over $\FF_5$.
\end{abstract}

\section{Introduction}\label{Sec:1}

Let $\FF_q$ denote the finite field of order $q$,
where $q$ is a prime power.
An $[n,k]_q$ code $C$ is a $k$-dimensional vector subspace
of $\FF_q^n$,
where $n$ and $k$ are called the length and the dimension of $C$,
respectively.
The weight $\wt(x)$ of a codeword $x$ is
the number of non-zero components of $x$.
The minimum non-zero weight of all codewords in $C$ is called
the minimum weight of $C$.
An $[n,k,d]_q$ code is an $[n,k]_q$ code with minimum weight $d$.

It is a fundamental problem  in coding theory to
determine the following values:
\begin{enumerate}
\item
the largest value $d_q(n,k)$ of $d$ for which there exists
an $[n,k,d]_q$ code.
\item
the smallest value $n_q(k,d)$ of $n$ for which there exists
an $[n,k,d]_q$ code.
\end{enumerate}
A code which achieves one of these two values is called 
{\em optimal}.
For $q \le 9$, the current knowledge on the 
values $d_q(n,k)$
can be obtained from~\cite{Grassl} 
(see also~\cite{Brouwer-Handbook} and \cite{mint}).
However, much work has been done concerning optimal codes 
for $q=2,3$ and $4$ only.
In this note, we consider optimal codes for $q=5$.
The smallest length $n$ for which  $d_5(n,k)$ is not determined
is $21$, more precisely,
$d_5(21,5)=13$ or $14$.

In this note, we demonstrate that there is no 
$[21,5,14]_5$ code.
The non-existence is established by
classifying codes with parameters
$[18,2,15]_5$ and $[18+t,2+t,14]_5$ ($t=0,1,2$).
The non-existence of a $[21,5,14]_5$ code
determines the following values. 

\begin{prop}\label{main}
$d_5(21+t,5+t)=13$ for $t=0,1,\ldots,4$.
\end{prop}

\begin{rem}
The above proposition yields that 
$n_5(5+t,14)=22+t$ for $t=0,1,\ldots,4$.
\end{rem}

The punctured code of an $[n,k,d]_5$ code with $d \ge 2$
is an $[n-1,k,d']_5$ code with $d' =d-1$ or $d$.
If there is an $[n,k,d]_5$ code then
there is an $[n-d,k-1,d']_5$ code with $d' \ge d/5$
(see~\cite[p.~302]{Brouwer-Handbook}).
Hence, as a consequence of the above proposition, we have
the following:

\begin{cor}
There is no code with parameters
\begin{align*}
&{[ 22+t, 5+t, 15 ]_5} \ (t=0,1,\dots,4), \\
&{[ 87+t, 6, 66+t ]_5} \ (t=0,1), \\
&{[ 88+t, 7, 66+t ]_5} \ (t=0,1), \\
&{[ 89+t, 8, 66+t ]_5} \ (t=0,1). 
\end{align*}
\end{cor}

Generator matrices of all codes given in this
note can be obtained electronically from

\medskip
\noindent
\url{http://yuki.cs.inf.shizuoka.ac.jp/codes/index.html}

\medskip
All computer calculations in this note were
done by programs in {\sc Magma}~\cite{Magma} and 
programs in the language C.

\section{Results}\label{Sec:Res}
\subsection{Method}

The covering radius of an $[n,k]_5$ code $C$
is the smallest integer $R$ such that spheres of radius $R$ 
around codewords of $C$ cover the space $\FF_5^n$.
A {\em shortened code} $C'$ of a code $C$ is the set of all codewords 
in $C$ which are $0$ in a fixed coordinate with that 
coordinate deleted.
A shortened code $C'$ of an $[n,k,d]_5$ code $C$ with $d \ge 2$
is an $[n-1,k,d]_5$ code if the deleted coordinate
is a zero coordinate and an $[n-1,k-1,d']_5$ 
code with $d' \ge d$ and covering radius $R\ge d-1$ 
otherwise.

Two $[n,k]_5$ codes $C$ and $C'$ are {\em equivalent} 
if there exists an $n \times n$ monomial matrix $P$ over $\FF_5$
with $C' = C \cdot P = \{ x P\:|\: x \in C\}$.  
To test equivalence of codes by a program in the language C, we use the
algorithm given in \cite[Section 7.3.3]{KO} as follows.
For an $[n,k]_5$ code $C$, define 
the digraph $\Gamma(C)$ with vertex set 
$C \cup (\{1,2,\dots,n\}\times (\FF_5-\{0\}))$
and arc set 
$\{(c,(j,c_j)), ((j,c_j),c) \mid c=(c_{1},\ldots,c_{n}) \in C,  1 \le j \le n\} 
\cup \{((j,y),(j,2y))\mid
1 \le j \le n,\ y \in \FF_5-\{0\}\}$.
Then, two $[n,k]_5$ codes $C$ and $C'$ are equivalent
if and only if $\Gamma(C)$ and $\Gamma(C')$  are isomorphic.
We use {\sc nauty}~\cite{nauty}
for digraph isomorphism testing.
It can be also done by the function {\tt IsIsomorphic} in {\sc Magma}
to test equivalence of codes.

An $[n,k,d]_5$ code $C$ gives $n$ shortened codes
and at least $k$ codes among them are $[n-1,k-1,d']_5$ codes
with $d' \ge d$.
Hence, 
by considering the inverse operation of shortening, 
any $[n,k,d]_5$ code with $d \ge 2$ is constructed from some
$[n-1,k-1,d']_5$ code with $d' \ge d$ and covering radius
$R\ge d-1$ as follows.
Let $C'$ be an $[n-1,k-1,d']_5$ code with $d' \ge d$.
Up to equivalence, 
we may assume that $C'$ has 
a generator matrix of the form $
\left(\begin{array}{cc}
I_{k-1} & A
\end{array}\right)$, where $I_{k-1}$ denotes the
identity matrix of order $k-1$.
Then, up to equivalence, 
an $[n,k,d]_5$ code, which is constructed from $C'$
by considering the inverse operation of shortening, 
has the following generator matrix
\begin{equation}\label{eq:S1}
\left(\begin{array}{ccc|c|cccccccc}
 &       &&0      & &  &\\
 &I_{k-1}&&\vdots & &A &\\
 &       &&0      & &  &\\
\hline
0& \cdots& 0 &1&b_1 &\cdots&b_{n-k}
\end{array}\right),
\end{equation}
where $b=(b_1,b_2,\ldots,b_{n-k}) \in\FF_5^{n-k}$
with $\wt(b) \ge d-1$.

\subsection{Non-existence of a $[21,5,14]_5$ code}

We remark that there is no code with parameters 
$[19+t,3+t,d \ge 15]_5$ $(t=0,1)$ and  $[18,2,d \ge 16]_5$ 
(see~\cite{Grassl}).
Thus, 
any $[19,3,14]_5$ code is constructed by (\ref{eq:S1}) from 
some $[18,2,14 \text{ or }15]_5$ code $C$ with covering radius
$R \ge 13$, and
any $[20+t,4+t,14]_5$ code is constructed 
by (\ref{eq:S1}) from 
some $[19+t,3+t,14]_5$ code $C$ with $R \ge 13$
$(t=0,1)$.

In order to determine whether there is a
$[21,5,14]_5$ code or not,
we classified codes
with parameters $[18,2,15]_5$ and $[18+t,2+t,14]_5$ $(t=0,1,2)$.
It is easy to see that
there is a unique $[18,2,15]_5$ codes, and 
there are ten $[18,2,14]_5$ codes, up to equivalence.
Using generator matrices in form (\ref{eq:S1})
of inequivalent $[18,2,14]_5$ codes and $[18,2,15]_5$ codes,
we constructed all $[19,3,14]_5$ 
codes which must be checked further for equivalences.
Similarly, from inequivalent $[19,3,14]_5$ codes,
we constructed all $[20,4,14]_5$ 
codes which must be checked further for equivalences.
By checking equivalences among these codes,
we completed a classification of $[19+t,3+t,14]_5$ codes $(t=0,1)$.

For the above parameters, 
the number $\#$ of inequivalent codes
is listed in Table~\ref{Tab:2}.
The number $\#_W$ of different weight enumerators and 
the number $\#_R$ of inequivalent codes with
covering radius $R$ are also listed.
Then we have the following:

\begin{prop}\label{2}
Every $[20,4,14]_5$ code has covering radius $12$ and
there is no $[21,5,14]_5$ code.
\end{prop}

Proposition \ref{2} completes the
proof of Proposition~\ref{main}.

\begin{table}[thb]
\caption{Non-existence of a $[21,5,14]_5$ code}
\label{Tab:2}
\begin{center}
{\small
\begin{tabular}{c|c|c|cc}
\noalign{\hrule height0.8pt}
Parameters & $\#$ &$\#_W$ & $\#_{\ge 13}$ & $\#_{12}$ \\
\hline
$[18,2,14]_5$ & 10  & 9& 10&0\\
$[18,2,15]_5$ &  1  & 1& 1&0\\
$[19,3,14]_5$ &572  &90& 572&0\\
$[20,4,14]_5$ &3564 & 727&0 &3564\\
$[21,5,14]_5$ &0    & --& --&--\\
\noalign{\hrule height0.8pt}
  \end{tabular}
}
\end{center}
\end{table}

\bigskip
\noindent
{\bf Acknowledgments.}
In this work, the supercomputer of ACCMS, Kyoto University
was partially used.
The authors would like to thank Markus Grassl
for useful comments~\cite{G}.
This work was supported by JST PRESTO program and 
JSPS KAKENHI Grant Number 23340021.



\begin{thebibliography}{99}
\bibitem{Magma}W. Bosma, J.J. Cannon, C. Fieker and A. Steel,
{\sl Handbook of Magma Functions (Edition 2.17)}, 2010, 5117 pages.


\bibitem{Brouwer-Handbook} A.E. Brouwer,
{``Bounds on the size of linear codes,''} {in Handbook of Coding Theory},
V.S. Pless and W.C. Huffman (Editors),
Elsevier, Amsterdam, 1998, pp.\ 295--461.


\bibitem{Grassl} M. Grassl,
Code tables: Bounds on the parameters of various types of codes,
Available online at
``\verb+http://www.codetables.de/+''.

\bibitem{G} M. Grassl, private communication, May 21, 2012.



\bibitem{KO}P. Kaski and  P.R.J. \"Osterg\aa rd,
{\sl Classification Algorithms for Codes and Designs}, 
Springer, Berlin, 2006.



\bibitem{nauty}B.D. McKay, 
nauty user's guide (version 2.4),
Available online at
``\verb+http://cs.anu.edu.au/people/bdm/nauty/+''.

\bibitem{mint}W.C. Schmid, and R. Sch\"urer,
MinT,
Available online at
``\verb+http://mint.sbg.ac.at/+''.

\end{thebibliography}
\end{document}